\documentclass{article}
\usepackage{amsmath}
\DeclareMathOperator{\Atang}{Atang.}
\begin{document}
\title{On the formation of continued fractions\footnote{Delivered to the St.--Petersburg Academy September 4, 1775.
Originally published as
{\em De formatione fractionum continuarum},
Acta Academiae Scientarum Imperialis Petropolitinae \textbf{3} (1782),
no. 1, 3--29, and
republished in \emph{Leonhard Euler, Opera Omnia}, Series 1:
Opera mathematica,
Volume 15, Birkh\"auser, 1992. A copy of the original text is available
electronically at the Euler Archive, at http://www.eulerarchive.org. This paper
is E522 in the Enestr\"om index.}}
\author{Leonhard Euler\footnote{Date of translation: August 12, 2005.
Translated from the Latin
by Jordan Bell, 3rd year undergraduate in Honours Mathematics, School of Mathematics and Statistics, Carleton University,
Ottawa, Ontario, Canada.
Email: jbell3@connect.carleton.ca.
This translation was written
during an NSERC USRA supervised by Dr. B. Stevens.
}}
\date{}
\maketitle

1. A universal principle for unfolding continued fractions is
found in the infinite series of quantities $A,B,C$, etc., of which
each third succeeds from the preceding two by a certain law,
either constant or variable inasmuch as they depend in turn upon
each other, so that it will be
\[
fA=gB+hC, f'B=g'C+h'D, f''C=g''D+h''E, f'''D=g'''E+h'''F, \textrm{ etc.}
\]
From here indeed are deduced the following equalities:
\begin{align*}
\frac{fA}{B}=&g+\frac{hC}{B}=g+\frac{f'h}{f'B:C}\\
\frac{f'B}{C}=&g'+\frac{h'D}{C}=g'+\frac{f''h'}{f''C:D}\\
\frac{f''C}{D}=&g''+\frac{h''E}{D}=g''+\frac{f'''h''}{f'''D:E}\\
\frac{f'''D}{E}=&g'''+\frac{h'''F}{E}=g'''+\frac{f''''h'''}{f''''E:F}\\
\textrm{etc.}&\qquad \textrm{etc.}
\end{align*}
But if now the latter values are substituted successively in place
of the prior ones, the following continued fraction
spontaneously emerges:
\[
\frac{fA}{B}=g+\frac{f'h}{g'+
\frac{f''h'}{g''+
\frac{f'''h''}{g'''+
\frac{f''''h'''}{g''''+\textrm{ etc.}}}}}
\]
whose value therefore is determined by the two first terms $A$ \& $B$ of
the series alone.

2. Therefore whenever such a progression of quantities $A,B,C,D,E$, etc.
is had, of which such a law is disposed, that each of the terms
depends in turn on the succeeding two by a certain law, then from this
a continued fraction may be deduced, whose value is able to be assigned.
Therefore if such a particular formula were disposed, such that the
expansion of it should lead to such a series of quantities $A,B,C,D,E$,
etc. of which each term is determined by the preceding two, 
from this continued fractions will be able to be derived, in what way
it will be able to be revealed most conveniently by several examples.

\begin{center}
{\Large I. The expansion of the formula\\
$s=x^n(\alpha-\beta x-\gamma xx)$}
\end{center}

3. In this formula the exponent $n$ is considered an indefinite,
successively receiving all the values $1,2,3,4,5,6$, etc., from which,
providing it will be $n>0$, this formula vanishes by putting $x=0$,
and then indeed it also vanishes by taking
\[
x=-\frac{\beta \pm \surd (\beta \beta+4\alpha \gamma)}{2\gamma}.
\]
With this having been noted, this formula will be differentiated so that it will become
\[
ds=n\alpha x^{n-1}dx-(n+1)\beta\int x^ndx-(n+2)\gamma \int x^{n+1}dx+s,
\]
from which by integrating the parts, this integration may then be
indicated as such 
\[
x=-\frac{\beta \pm \surd (\beta \beta+4\alpha \gamma)}{2\gamma}.
\]
Then, if after this integration, having been completed
such that the integral should vanish by putting $x=0$, it should be
set
\[
n\alpha \int x^{n-1}dx=(n+1)\beta \int x^n dx+(n+2)\gamma \int x^{n+1}dx,
\]
of course in which case it would be $s=0$, it will be
\[
n\alpha \int x^{n-1}dx=(n+1)\beta \int x^ndx+(n+2)\gamma \int x^{n+1}dx,
\]
which is a relation of this type between three integral formulas
following each other, which we desire for the formation of
a continued fraction, seeing that these integral formulas, if
in place of $n$ is written successively the number $1,2,3,4,5,6$, etc.,
the quantities $A,B,C,D$, etc. are given to us.

4. We write therefore in place of $n$ the consecutive numbers
$1,2,3,4$, etc., so that these relations may be produced:
\begin{align*}
\alpha \int dx=2\beta \int xdx+3\gamma \int xxdx\\
2\alpha \int xdx=3\beta \int xxdx+4\gamma \int x^3dx\\
3\alpha \int xxdx=4\beta \int x^3dx+3\gamma \textrm{[sic]} \int x^4dx\\
4\alpha \int x^3dx=5\beta \int x^4dx+6\gamma \int x^5dx\\
\qquad \textrm{ etc.} \qquad \textrm{ etc.}
\end{align*}
Then we will therefore have
\begin{align*}
A=\int dx=x=-\frac{\beta \pm \surd (\beta \beta+4\alpha \gamma)}{2\gamma},\\
B=\int xdx=\frac{1}{2}xx=\frac{1}{2}\Big(\frac{-\beta \pm \surd(\beta \beta+4\alpha \gamma)}{2\gamma}\Big)^2,\\
C=\int xxdx=\frac{1}{3}x^3, \quad D=\int x^3dx=\frac{1}{4}x^4\\
\qquad \textrm{ etc.} \qquad \textrm{ etc.}
\end{align*}
Thereupon indeed, for the letters $f,g,h$ will be had these values:
\begin{align*}
f=\alpha, \, f'=2\alpha, \, f''=3\alpha, \, f'''=4\alpha \textrm{ etc.}\\
g=2\beta, \, g'=3\beta, \, g''=4\beta, \, g'''=5\beta \textrm{ etc.}\\
h=3\gamma, \, h'=4\gamma, \, h''=5\gamma, \, h'''=6\gamma \textrm{ etc.},
\end{align*}
from which values the following continued fraction results:
\[
\frac{\alpha A}{B}=2\beta+
\frac{6\alpha \gamma}{3\beta+
\frac{12\alpha \gamma}{4\beta+
\frac{20\alpha \gamma}{5\beta+
\frac{30\alpha \gamma}{6\beta+\textrm{ etc.}}}}},
\]
whose value therefore is
\[
\frac{4\alpha \gamma}{-\beta + \surd(\beta \beta+4\alpha \gamma)}=
\beta+\surd(\beta \beta +4\alpha \gamma).
\]

5. So that this continued fraction should be rendered more elegant,
in place of $\alpha \gamma$ we write $\frac{1}{2}\delta$, and it will
proceed
\[
\beta+\surd(\beta \beta+2\delta)=2\beta+
\frac{3\delta}{3\beta+
\frac{6\delta}{4\beta+
\frac{10\delta}{5\beta+
\frac{15\delta}{6\beta+\textrm{ etc.}}}}}.
\]
Since moreover the head of this expression is seen to have been truncated,
with this head having been added we may put
\[
s=\beta+\frac{\delta}{2\beta+
\frac{3\delta}{3\beta+
\frac{6\delta}{4\beta+
\frac{10\delta}{5\beta+\textrm{ etc.}}}}},
\]
and it will be
\[
s=\beta+\frac{\delta}{\beta+\surd(\beta \beta+2\delta)},
\]
which expression is reduced to this:
\[
s=\frac{1}{2}\beta+\frac{1}{2}\surd (\beta \beta+2\delta).
\]

6. This continued fraction is in fact able to be reduced to even
greater simplicity, if in place of $\delta$ we write $2\varepsilon$, so that
it would be
\[
\frac{1}{2}\beta+\frac{1}{2}\surd(\beta \beta+4\varepsilon)=\beta+
\frac{2\varepsilon}{2\beta+
\frac{6\varepsilon}{3\beta+
\frac{12\varepsilon}{4\beta+
\frac{20\varepsilon}{5\beta+20 \textrm{[sic]  etc.}}}}}.
\]
And if now the first fraction is depressed by 2,
the second by 3, the third by 4, the fourth by 5, etc., the following
form will be produced:
\[
\frac{1}{2}\beta+\frac{1}{2}\surd(\beta \beta+4\varepsilon)=\beta+
\frac{\varepsilon}{\beta+
\frac{\varepsilon}{\beta+
\frac{\varepsilon}{\beta+
\frac{\varepsilon}{\beta+\textrm{ etc.}}}}},
\]
which is of the greatest simplicity, and if its sum is considered an unknown,
and it is denoted as equal to $z$,
it will be in turn be $z=\beta+\frac{\varepsilon}{z}$,
and thus $zz=\beta z+\varepsilon$, from which it would be
$z=\frac{\beta+\surd(\beta \beta+4\varepsilon)}{2}$, which is the same.

7. Indeed this most simple sum is able to be immediately
deduced from the formula taken initially
\[
s=x^n(\alpha-\eta x-\gamma xx),
\]
and since we have put this equal to nothing, it will
certainly be
$\alpha=\beta x+\gamma xx$, and in the very same way
\[
\alpha x=\beta xx+\gamma x^3, \, \alpha xx=\beta x^3+\gamma x^4, \textrm{ etc.},
\]
thus so that for the series $A,B,C,D$, etc. we may have this
simple series of powers: $1,x,x^3,x^3,x^4$, etc.,
then indeed each of the letters
$f,g,h$, etc. becomes $\alpha,\beta,\gamma$, etc., from which
arises this continued fraction:
\[
\frac{\alpha}{x}=\beta+\frac{\alpha \gamma}{\beta+
\frac{\alpha \gamma}{\beta+
\frac{\alpha \gamma}{\beta+\textrm{ etc.}}}},
\]
where it is $\frac{1}{x}=\frac{\beta+\surd(\beta \beta+4\alpha \gamma)}{2\alpha}$.
Therefore the value of this fraction is $\frac{1}{2}\beta+
\frac{1}{2}\surd(\beta \beta+4\alpha \gamma)$, such as before,
since $\alpha \gamma=\varepsilon$.

\begin{center}
{\Large II. The expansion of the formula\\
$s=x^n(a-x)$}
\end{center}

8. This formula will therefore vanish by setting $x=a$.
Moreover it is $ds=nax^{n-1}dx-(n+1)x^ndx$, which expression
is comprised by two terms; it may be reduced to a fraction,
whose denominator is $\alpha+\beta x$, so that it will become
\[
ds=\frac{na\alpha x^{n-1}dx+(\beta na-\alpha(n+1))x^n dx-\beta (n+1)x^{n+1}
dx}{\alpha+\beta x}.
\]
Therefore by integrating each of the members, it will be
\[
s=na\alpha \int \frac{x^{n-1}dx}{\alpha+\beta x}+
(n\beta a-(n+1)\alpha) \int \frac{x^n dx}{\alpha+\beta x}-
\beta (n+1)\int \frac{x^{n+1} dx}{\alpha+\beta x},
\]
which if after each are integrated we set $x=a$, so that
it may become $s=0$, we will have this reduction:
\[
na\alpha \int \frac{x^{n-1}dx}{\alpha+\beta x}=((n+1)\alpha-n\beta a)
\int \frac{x^n dx}{\alpha+\beta x}+(n+1)\beta \int \frac{x^{n+1}dx}{\alpha+\beta x}).
\]

9. In place of $n$ we shall now successively substitute the numbers
$1,2,3,4$, etc., and then by comparison with the general formula
that has been established we will have
\[
A=\int \frac{dx}{\alpha+\beta x}, \qquad B=\int \frac{xdx}{\alpha+\beta x}, \qquad C=\int \frac{xxdx}{\alpha+\beta x}, \textrm{ etc.}
\]
where indeed after integration it ought to be made $x=a$. 
Thereafter truly we will have
\begin{align*}
f=a\alpha, \, f'=2a\alpha, \, f''=3a\alpha, \, f'''=4a\alpha, \textrm{ etc.}\\
g=2\alpha-\beta a, \, g'=3\alpha-2\beta a, \, g''=4\alpha-3\beta a, \, \textrm{ etc.}\\
h=2\beta, h'=3\beta, \, h''=4\beta, \, h'''=5\beta, \textrm{ etc.}
\end{align*}
and thus from this arises the following continued fraction:
\[
\frac{\alpha aA}{B}=(2\alpha-\beta a)+\frac{4a\alpha \beta}{(3\alpha-2\beta a)+
\frac{ga\alpha \beta}{(4\alpha-3\beta a)+
\frac{16a\alpha \beta}{(5\alpha-4\beta a)+\textrm{ etc.}}}}.
\]

10. By integrating, it may moreover be established
\[
\int \frac{dx}{\alpha+\beta x}=\frac{1}{\beta}l\frac{\alpha+\beta x}{\alpha},
\]seeing that the integral ought to vanish by making it $x=0$.
Now therefore it may be $x=a$, and it will hence be $A=\frac{1}{\beta}l\frac{\alpha+\beta a}{\alpha}$. In turn
\[
\int \frac{xdx}{\alpha+\beta x}=\frac{1}{\beta}\Big(x-\frac{\alpha}{\beta}l\frac{\alpha+\beta x}{\alpha}\Big),
\]
and by it being made $x=a$ it shall become
\[
B=\frac{a}{\beta}-\frac{\alpha}{\beta \beta}l\frac{\alpha+\beta a}{\alpha},
\]
wherefore the value of our continued fraction will be
\[
\frac{\alpha a\beta l\frac{\alpha+\beta a}{\alpha}}{a\beta-\alpha l\frac{\alpha+\beta a}{\alpha}};
\]
moreover it is evident for nothing of universality to perish, even if
it is set $a=1$; then in fact it will be
\[
\frac{\alpha \beta l\frac{\alpha+\beta}{\alpha}}{\beta-\alpha l\frac{\alpha+\beta}{\alpha}}=(2\alpha-\beta)+\frac{4\alpha \beta}{(3\alpha-2\beta)+
\frac{9\alpha \beta}{(4\alpha-3\beta)+\textrm{ etc.}}}.
\]

11. Moreover, the whole of this expression manifestly depends singularly
on the ratio of the numbers $\alpha$ and $\beta$; from this, we may
take $\alpha=1$ and $\beta=n$, and then this continued fraction will be
seen:
\[
\frac{nl(1+n)}{n-l(1+n)}=(2-n)+\frac{4n}{(3-2n)+
\frac{9n}{(4-3n)+
\frac{16n}{(5-4n)+\textrm{ etc.}}}},
\]
for which if we set this series after the law $1+n$ and
set the sum equal to $s$, so that it will be
\[
s=1+\frac{n}{(2-n)+
\frac{4n}{(3-2n)+
\frac{9n}{(4-3n)+
\frac{16n}{(5-4n)+\textrm{ etc.}}}}},
\]
it will be
\[
s=1+\frac{n(n-l(1+n))}{nl(1+n)}=1+\frac{n-l(1+n)}{l(1+n)}=
\frac{n}{l(1+n)}.
\]

12. We shall run through several examples, and the first shall be
$n=1$, where it will be
\[
\frac{1}{l2}=1+
\frac{1}{1+
\frac{4}{1+
\frac{9}{1+
\frac{16}{1+\textrm{ etc.}}}}}.
\]
On the other hand by putting $n=2$ it will be
\[
\frac{2}{l3}=1+
\frac{2}{0+
\frac{8}{-1+
\frac{18}{-2+
\frac{32}{-3+
\frac{50}{-4+\textrm{ etc.}}}}}},
\]
which expression however, on account of the negative quantities,
is not very pleasant; insofar as it will turn out that whenever $n>1$,
it will be more worth the work to unfold these cases, than when $n$ is
taken as less than unity.

13. Therefore it is able to be made easily;
 it may be reverted to an expression containing the letters
$\alpha$ and $\beta$, and then by supplying the head because it is missing,
this form will be produced:
\[
\frac{\beta}{l\frac{\alpha+\beta}{\alpha}}=\alpha+
\frac{\alpha \beta}{(2\alpha-\beta)+
\frac{4\alpha \beta}{(3\alpha-2\beta)+
\frac{9\alpha \beta}{(4\alpha-3\beta)+\textrm{ etc.}}}}. 
\]
We now put $\alpha=n-m$ and $\beta=2m$,\footnote{Translator: Original reads ``$n=n-m$''.}
so that we may obtain
the following form:
\[
\frac{2m}{l\frac{n+m}{n-m}}=
n-m+\frac{2m(n-m)}{2n-4m+
\frac{8m(n-m)}{3n-7m+
\frac{18m(n-m)}{4n-10m+\textrm{ etc.}}}},
\]
from which the following special cases are deduced.

If $m=1$ and $n=3$ it will be
\[
\frac{2}{l2}=2+\frac{4}{2+
\frac{16}{2+
\frac{36}{2+
\frac{64}{2+\textrm{ etc.}}}}},
\]
which fraction divided by 2 and reduced is rendered as such:
\[
\frac{1}{l2}=1+\frac{1}{1+
\frac{4}{1+
\frac{9}{1+
\frac{16}{1+\textrm{ etc.}}}}},
\]
which here is arrived at like before.

Were it $m=1$ and $n=4$ it will be
\begin{align*}
\frac{2}{l\frac{5}{2}}=3+\frac{6}{4+
\frac{24}{5+
\frac{54}{6+
\frac{96}{7+\textrm{ etc.}}}}}\\
=
3+\frac{6\cdot 1}{4+
\frac{6\cdot 4}{5+
\frac{6\cdot 9}{6+
\frac{6\cdot 16}{7+\textrm{ etc.}}}}}.
\end{align*}
Were it $m=1$ and $n=5$, it will be
\[
\frac{2}{l\frac{3}{2}}=4+\frac{8}{6+
\frac{32}{8+
\frac{72}{10+
\frac{128}{12+\textrm{ etc.}}}}},
\]
or
\begin{align*}
\frac{1}{l\frac{3}{2}}=2+\frac{2}{3+
\frac{8}{4+
\frac{18}{5+
\frac{32}{6+\textrm{ etc.}}}}}\\
=2+\frac{2\cdot 1}{3+
\frac{2\cdot 4}{4+
\frac{2\cdot 9}{5+
\frac{2\cdot 16}{6+\textrm{ etc.}}}}}.
\end{align*}

\begin{center}
{\Large III. The expansion of the formula\\
$s=x^n(1-x^2)$}
\end{center}

14. This formula therefore vanishes in the cases $x=0$ and $x=1$.
Then indeed, it shall be
\[
ds=nx^{n-1}dx-(n+2)x^{n+1}dx,
\]
which differential uplifted with the denominator $\alpha+\beta xx$ would
be
\[
ds=\frac{n\alpha x^{n-1}dx+(n\beta-(n+2)\alpha)x^{n+2}dx-
(n+2)\beta x^{n+3}dx}{\alpha+\beta xx};
\]
here, by now integrating this in turn it will be
\[
s=n\alpha \int \frac{x^{n-1}dx}{\alpha+\beta xx}+(n\beta-(n+2)\alpha)
\int \frac{x^{n+1}dx}{\alpha+\beta xx}-(n+2)\beta \int \frac{x^{n+3}dx}{\alpha+\beta xx}.
\]
And if now after these integrations it is set $x=1$, this reduced integral
will be produced:
\[
n\alpha \int \frac{x^{n-1}dx}{\alpha+\beta xx}=((n+2)\alpha-n\beta)
\int \frac{x^{n+1}dx}{\alpha+\beta xx}+(n+2)\beta \int \frac{x^{n+3}dx}{\alpha+\beta xx}.
\]

15. Because these powers of $x$ are augmented two by two,
for the exponent $n$ we shall assign successively the values
$1,3,5,7,9$, etc., and it shall be set:
\[
A=\int \frac{dx}{\alpha+\beta xx}, \qquad
B=\int \frac{xxdx}{\alpha+\beta xx}, \qquad
C=\int \frac{x^4dx}{\alpha+\beta xx}, \textrm{ etc.}
\]
Then indeed from these, the letters $f,g,h$ will be derived:
\begin{align*}
f=\alpha, \, f'=3\alpha, \, f''=5\alpha, \, f'''=7\alpha, \textrm{ etc.}\\
g=3\alpha-\beta, \, g'=5\alpha-3\beta, \, g''=7\alpha-5\beta, \textrm{ etc.}\\
h=3\beta, h'=5\beta, \, h''=7\beta, \, h'''=9\beta, \textrm{ etc.},
\end{align*}
from which is born the following continued fraction:
\[
\frac{\alpha A}{B}=3\alpha-\beta+\frac{9\alpha \beta}{5\alpha-3\beta+
\frac{25\alpha \beta}{7\alpha-5\beta+
\frac{49\alpha \beta}{9\alpha-7\beta+\textrm{ etc.}}}}.
\]

16. Because it is $B=\int \frac{xxdx}{\alpha+\beta xx}$, it will be
$B=\frac{1}{\beta} \int dx- \frac{\alpha}{} \textrm{ [sic] } \int \frac{dx}{\alpha+\beta xx}$,
and thus $B=\frac{1}{\beta}-\frac{\alpha}{\beta}A$,
by whose value being substituted in we will have
\[
\frac{\alpha \beta A}{1-\alpha A}=3\alpha-\beta+\frac{9\alpha \beta}{5\alpha-3\beta+
\frac{25\alpha \beta}{7\alpha-5\beta+\textrm{ etc.}}},
\]
which, because it is missing its head, we will set in front $\alpha+\beta+
\alpha \beta$; then moreover the sum will be $\beta+\frac{1}{A}$,
so that we will thus have
\[
\beta+\frac{1}{A}=\alpha+\beta+\frac{\alpha \beta}{3\alpha-\beta+
\frac{9\alpha \beta}{5\alpha-3\beta+
\frac{25\alpha \beta}{7\alpha-5\beta+\textrm{ etc.}}}}.
\]
with it arising that $A=\int \frac{dx}{\alpha+\beta xx}$,
with this integral being taken so that it vanishes by putting $x=0$,
and indeed also by making it $x=1$.

17. First we shall expand the simplest case, in which $\alpha=1$ and
$\beta=1$, where it will be $A=\frac{\pi}{4}$, for
which we will have
\[
1+\frac{4}{\pi}=2+\frac{1}{2+
\frac{9}{2+
\frac{25}{2+
\frac{49}{2+\textrm{ etc.}}}}},
\]
that is it will be
\[
\frac{4}{\pi}=1+\frac{1}{2+
\frac{9}{2+
\frac{25}{2+\textrm{ etc.}}}},
\]
which is the same continued fraction produced first before by Brouncker,
the investigation of which, which was elicited by Wallis through greatly
tedious calculations, here proceeds immediately from our formula by itself.

18. Our general form as well provides infinite other similar
expressions, precisely as the letters $\alpha$ and $\beta$ are taken
in varied ways. And firstly indeed, if $\alpha$ and $\beta$
were positive numbers, the value of the letter $A$ is always able
to be expressed by circular arcs, and conversely indeed by logarithms.
Were it therefore $\beta=1$, it will be
\[
A=\int \frac{dx}{\alpha+xx}=\frac{1}{\surd \alpha} \Atang \frac{x}{\surd \alpha}=
\frac{1}{\surd \alpha} \Atang \frac{1}{\surd \alpha},
\] 
from which is born this continued fraction:
\[
1+\frac{\surd \alpha}{\Atang \frac{1}{\surd \alpha}}=\alpha+1+
\frac{\alpha}{3\alpha-1+
\frac{9\alpha}{5\alpha-3+
\frac{25\alpha}{7\alpha-5+\textrm{ etc.}}}}.
\]
Here therefore it should be taken $\alpha=3$, and because $\Atang \frac{1}{\surd 3}=
\frac{\pi}{3}$ we will have
\[
1+\frac{6\surd 3}{\pi}=4+\frac{3}{8+
\frac{27}{12+
\frac{75}{16+
\frac{147}{20+\textrm{ etc.}}}}},
\]
or
\[
1+\frac{6\surd 3}{\pi}=4+\frac{3\cdot 1}{8+
\frac{3\cdot 9}{12+
\frac{3\cdot 25}{16+
\frac{3\cdot 49}{20+\textrm{ etc.}}}}}.
\]

19. Now $B$ shall be a particular positive number, and indeed it is
\[
A=\int \frac{dx}{\alpha+\beta xx}=\frac{1}{\beta}\int \frac{dx}{\frac{\alpha}{\beta}+xx},
\]
which having been integrated will be $A=\frac{1}{\surd \alpha \beta}\Atang \surd \frac{\beta}{\alpha}$. Then therefore we will have
\[
\beta+\frac{\surd \alpha \beta}{\Atang \surd \frac{\beta}{\alpha}}=\alpha+\beta+
\frac{\alpha \beta}{3\alpha-\beta+
\frac{9\alpha \beta}{5\alpha-3\beta+\textrm{ etc.}}}.
\]
We shall now therefore set $\alpha+\beta=2n$ and $\alpha-\beta=2m$, so that
it shall be $\alpha=m+n$ and $\beta=n-m$, with whose values having been
put in it will be
\[
n-m+\frac{\surd(nn-mm)}{\Atang \surd \frac{n-m}{n+m}}=2n+\frac{nn-mm}{2n+4m+
\frac{9(nn-mm)}{2n+8m+\textrm{ etc.}}}.
\]

20. We may consider as well the case, in which $\beta$ is a negative
number, and by putting $\beta=-\gamma$ it will be
\[
A=\int \frac{dx}{\alpha-\gamma xx}=\frac{1}{\surd} \int
\frac{dx}{\frac{\alpha}{\surd}-xx} \textrm{[sic]},
\]
whose integral is
\[
A=\frac{1}{2\surd \alpha \gamma}l\frac{\surd{\frac{\alpha}{\gamma}}+x}
{\surd{\frac{\alpha}{\gamma}}-x};
\]
with it having been made therefore $x=1$ it will be
\[
A=\frac{1}{2\surd \alpha \gamma}l \frac{\surd \alpha+\surd \gamma}{\surd \alpha -\surd \gamma},
\]
from which is born this continued fraction:
\[
-\gamma+\frac{2\surd \alpha \gamma}{l \frac{\surd \alpha+\surd \gamma}{\surd \alpha-\surd \gamma}}=\alpha-\gamma-\frac{\alpha \gamma}{3\alpha+\gamma-
\frac{9\alpha \gamma}{5\alpha+3\gamma-
\frac{25\alpha \gamma}{7\alpha+5\gamma-\textrm{ etc.}}}},
\]
and thus we have made in this way new continued fractions, of which
the values may moreover be exhibited by logarithms,
and which are altogether different from that which we came upon
before.

21. This case presents itself being more worthy of notice than
the others,
when $\gamma=\alpha$. Or, because it turns out the same, $\alpha=1$
and $\gamma=1$; because then it is $l\frac{\surd \alpha + \surd
\gamma}{\surd \alpha- \surd \gamma}=l\infty=\infty$, we will have
\[
-1=0-\frac{1}{4-
\frac{9}{8-
\frac{25}{12-\textrm{ etc.}}}},
\] 
or by changing the sign
\[
1=\frac{1}{4-
\frac{9}{8-
\frac{25}{12-
\textrm{ etc.}}}},
\]
whence the first denominator
\[
4-\frac{9}{8-
\frac{25}{12-\textrm{ etc.}}}
\]
ought to be equal to 1. Therefore it will be
\[
0=3-\frac{9}{8-
\frac{25}{12-\textrm{ etc.}}},
\]
or
\[
1=\frac{3}{8-
\frac{25}{12-\textrm{ etc.}}},
\]
where the denominator ought to be equal to 3, from which it shall be
\[
0=5-\frac{25}{12-\textrm{ etc.}},
\]
whose denominator ought to be equal to 5, from which it shall be
\[
0=7-\frac{49}{16-
\frac{81}{20-\textrm{ etc.}}},
\]
from which the true order is easily perceived.

22. We take $\alpha=4$ and $\gamma=1$, and this fraction will be
born
\[
-1+\frac{4}{l3}=3-\frac{4\cdot 1}{13-
\frac{4\cdot 9}{23-
\frac{4\cdot 25}{33-
\frac{4\cdot 49}{43-\textrm{ etc.}}}}}
\]
But if on the other hand we let $\alpha=9$ and $\gamma=1$ it will
be
\[
-1+\frac{6}{l2}=8-\frac{9\cdot 1}{28-
\frac{9\cdot 9}{48-
\frac{9\cdot 25}{68-
\frac{9\cdot 49}{88-\textrm{ etc.}}}}}.
\]

\begin{center}
{\Large IV. The expansion of the formula\\
$s=x^ne^{\alpha x}(1-x)$}
\end{center}

22.[sic] This $e$ denotes the number whose hyperbolic logarithm is unity, thus
so that $d.e^{\alpha x}=\alpha dx e^{\alpha x}$. Then it will
therefore be
\[
ds=nx^{n-1}dxe^{\alpha x}+(\alpha-(n+1))x^ndxe^{\alpha x}-
\alpha x^{n+1}dxe^{\alpha x},
\]
from which in turn by integrating it will be
\[
s=n\int x^{n-1}dxe^{\alpha x}+(\alpha-(n+1))\int x^ndxe^{\alpha x}-
\alpha \int x^{n+1}dxe^{\alpha x}.
\]
But if now after integration it is put $x=1$, it will be
\[
n\int x^{n-1} dxe^{\alpha x}=(n+1-\alpha)\int x^ndxe^{\alpha x}+
\alpha \int x^{n+1}dxe^{\alpha x}.
\]

23. But if now in place of $n$ we successively write the numbers
$1,2,3,4$, and so on, we may make
\begin{align*}
A\int e^{\alpha x}dx=\frac{1}{\alpha}(e^{\alpha}-1) \, \textrm{ and } \,
B=\int xdxe^{\alpha x}=\frac{\alpha-1}{\alpha \alpha}e^{\alpha}+
\frac{1}{\alpha \alpha}\\
f=1, \, f'=2, \, f''=3, \, f'''=4, \textrm{ etc.}\\
g=2-\alpha, \, g'=3-\alpha, \, g''=4-\alpha, \textrm{ etc.}\\
h=\alpha, \, h'=\alpha, \, h''=\alpha, \, h'''=\alpha, \textrm{ etc.}
\end{align*}
and this continued fraction will be produced:
\[
\frac{A}{B}=2-\alpha+\frac{2\alpha}{3-\alpha+
\frac{3\alpha}{4-\alpha+
\frac{4\alpha}{5-\alpha+\textrm{ etc.}}}}.
\]
We now adjoin at the head $1-\alpha+\alpha$, whose value will be
\[
1-\alpha+\frac{(\alpha-1)e^\alpha+1}{e^\alpha-1}=\frac{\alpha}{e^\alpha-1},
\]
from which is obtained this quite elegant continued fraction:
\[
\frac{\alpha}{e^\alpha-1}=1-\alpha+\frac{\alpha}{2-\alpha+
\frac{2\alpha}{3-\alpha+
\frac{3\alpha}{4-\alpha+\textrm{ etc.}}}},
\]
from which it is apparent that if it were $\alpha=0$, because
$e^\alpha-1=\alpha$, for it to be certainly $1=1$.

24. We may consider several particular cases; 
and first, if it were $\alpha=1$, it will be
\[
\frac{1}{e-1}=0+\frac{1}{1+
\frac{2}{2+
\frac{3}{3+
\frac{4}{4+\textrm{ etc.}}}}},
\]
which continued fraction is easily transformed into this:
\[
\frac{1}{e-1}=\frac{1}{1+
\frac{\frac{1}{1}}{1+
\frac{\frac{1}{2}}{1+
\frac{\frac{1}{3}}{1+
\frac{\frac{1}{4}}{1+\textrm{ etc.}}}}}},
\]
from which it is
\[
e-1=1+\frac{\frac{1}{1}}{1+
\frac{\frac{1}{2}}{1+
\frac{\frac{1}{3}}{1+\textrm{ etc.}}}}.
\]
Moreover, this in turn having been resolved by partial fractions gives
\[
e-1=+\frac{1}{1+
\frac{1}{2+
\frac{2}{3+
\frac{3}{4+
\frac{4}{5+\textrm{ etc.}}}}}},
\]
from which follows
\[
\frac{1}{e-2}=1+\frac{1}{2+
\frac{2}{3+
\frac{3}{4+
\frac{4}{4 \textrm{ [sic] }+\textrm{etc.}}}}},
\]
which forms because of their great simplicity are noteworthy. From the
second last, it may be
\[
e=2+\frac{1}{1+
\frac{1}{2+
\frac{2}{3+
\frac{3}{4+\textrm{ etc.}}}}};
\]
by taking successively $1,2,3$, and so on for the members, the following
approximations arise:
\begin{align*}
e=2,0000\\
e=3,0000\\
e=2,6666\\
e=2,7272\\
e=2,7169
\end{align*}
which values, alternately greater and smaller, converge to the truth
readily enough.

25. We may take $\alpha=2$ and it will be
\[
\frac{2}{ee-1}=-1+\frac{2}{0+
\frac{4}{1+
\frac{6}{2+
\frac{8}{3+\textrm{ etc.}}}}}.
\]
From this fraction in turn is deduced this:
\[
\frac{2(ee-1)}{ee+1}=0+\frac{4}{1+
\frac{6}{2+
\frac{8}{3+\textrm{ etc.}}}},
\]
and in a similar way, if larger numbers are taken for $\alpha$,
a reduction will be able to be made.

26. As well, negative numbers are able to be taken for $\alpha$.
In this way, if it were $\alpha=-1$, it will become
\[
\frac{e}{e-1}=2-\frac{1}{3-
\frac{2}{4-
\frac{3}{5-
\frac{4}{6-\textrm{ etc.}}}}},
\]
which is reduced to this form:
\[
\frac{e}{e-1}=2+\frac{1}{-3+
\frac{2}{4+
\frac{3}{-5+
\frac{4}{6+\textrm{ etc.}}}}},
\]
and in a similar way larger values are able to be dealt with.

27. Now too we set $\alpha=\frac{1}{2}$, and
this expression will be found:
\[
\frac{1}{2(\surd e-1)}=\frac{1}{2}+\frac{\frac{1}{2}}{\frac{3}{2}+
\frac{1}{\frac{5}{2}+
\frac{\frac{3}{2}}{\frac{7}{2}+
\frac{\frac{5}{2}}{\frac{9}{2}+\textrm{ etc.}}}}}
\]
which reduced by partial fractions comes out as
\[
\frac{1}{-1+\surd e}=1+\frac{2}{3+
\frac{4}{5+
\frac{6}{7+
\frac{8}{9+\textrm{ etc.}}}}}.
\]
In a similar way if we take $\alpha=\frac{1}{3}$, it will be
\[
\frac{1}{3(\sqrt[3] e-1)}=2:3+\frac{1:3}{5:3+
\frac{2:3}{8:3+
\frac{3:3}{11:3+
\frac{4:3}{14:3+\textrm{ etc.}}}}},
\]
which reduced by partial fractions gives
\[
\frac{1}{-1+\sqrt[3] e}=2+\frac{3}{5+
\frac{6}{8+
\frac{9}{11+
\frac{12}{14+\textrm{ etc.}}}}}.
\]
And if it is put $\alpha=\frac{2}{3}$, this continued fraction will
be produced:
\[
\frac{2}{3(\sqrt[3] (ee-1)}=1:3+
\frac{2:3}{4:3+
\frac{4:3}{7:3+
\frac{6:3}{10:3+
\frac{8:3}{13:3+\textrm{ etc.}}}}},
\]
which reduced by partial fractions will be
\[
\frac{2}{\sqrt[3] (ee-1)}=1+\frac{6}{4+
\frac{12}{7+
\frac{18}{10+
\frac{24}{13+\textrm{ etc.}}}}}.
\]

28. With these fairly simple principles disclosed, in a similar
way it will be 
permitted to treat other more general ones, which will be led to
continued fractions much more abstrusely, such that this will be accessible
from the cases which follow. 

\begin{center}
{\Large V. The expansion of the formula\\
$s=x^n(a-bx^\theta-cx^{2\theta})^\lambda$}
\end{center} 

29. Here therefore it will be
\[
ds=(a-bx^\theta-cx^{2\theta})^{\lambda-1}(nax^{n-1}dx-b(n+\lambda \theta)
x^{n+\theta-1}dx-c(n+2\lambda \theta)x^{n+2\theta-1}dx,
\]
from which by its parts being integrated, then indeed by putting $a-bx^\theta-
cx^{2\theta}=0$, (insofar as if it were it will be $x^\theta=
-\frac{b+\surd(bb+4ac)}{2c}$) this reduction will be had in general:
\begin{align*}
na\int x^{n-1}dx(a-bx^\theta-cx^{2\theta})^{\lambda-1}\\
=(n+\lambda \theta)b\int x^{n+\theta-1}dx(a-bx^\theta-cx^{2\theta})^{\lambda-1}\\
+(n+2\lambda \theta)c\int x^{n+2\theta-1}dx(a-bx^\theta-cx^{2\theta})^{\lambda-1}.
\end{align*} 

30. But if we want to compare this form with our general one related
earlier, we must discover the values for the letter
$n$ which are to be successively assumed for different $\theta$.
Then it is not necessary that the first value for $n$,
as we have so far made it, be taken as equal to 1; we shall set
therefore the first value of it as equal to $\alpha$, and we shall
search for the values of the two following integral formulas, namely:
\begin{align*}
A=\int x^{\alpha-1}dx(a-bx^\theta-cx^{2\theta})^{\lambda-1} \textrm{ and}\\
B=\int x^{\alpha+\theta-1}dx(a-bx^\theta-cx^{2\theta})^{\lambda-1},
\end{align*}
which integrals are taken such that they vanish by putting
$x=0$, which having been made the former value of $x$ ought to be taken,
which returns the formula $a-bx^\theta-c\alpha^{2\theta}=0$. Since however
it is not permitted for this to be carried out in general,
we shall be content to indicate these values by the letters $A$ and $B$,
which therefore consider as known.

31. Thereafter indeed the letters $f,g,h$, with this having been derived
will take on the following values:
\begin{align*}
f=\alpha a, \, f'=(\alpha+\theta)a, \, f''=(\alpha+2\theta)a, \, f'''=(\alpha+3\theta)a, \textrm{ etc.}\\
g=(\alpha+\lambda \theta)b, \, g'=(\alpha+\theta+\lambda \theta)b, \, g''=(\alpha+2\theta+\lambda \theta)b, \textrm{ etc.}\\
h=(\alpha+2\lambda \theta)c, \, h'=(\alpha+\theta+2\lambda \theta)c, \, h''=(\alpha+2\theta+2\lambda \theta)c, \textrm{ etc.}
\end{align*}
From this therefore will be formed the following continued fraction:
\[
\frac{\alpha aA}{B}=(\alpha+\lambda \theta)b+\frac{(\alpha+\theta)(\alpha+\lambda \theta)ac}{(\alpha+\theta+\lambda \theta)b+
\frac{(\alpha+2\theta)(\alpha+\theta+2\lambda \theta)ac}{(\alpha+2\theta+\lambda \theta)b+
\frac{(\alpha+3\theta)(\alpha+2\theta+2\lambda \theta)ac}{(\alpha+3\theta+
\lambda \theta)b \textrm{ etc.}}}},
\]
which form is certainly most general, of which however we will not
exclude another
derivation.  

\begin{center}
{\Large VI. The expansion of the formula\\
$s=x^n(1-x^\theta)^\lambda$}
\end{center}

32. Here therefore it would be
\[
ds=nx^{n-1}dx(1-x^\theta)^\lambda-\lambda \theta x^{n+\theta-1}dx(1-x^\theta)^{\lambda-1},
\] 
from which
will emerge two integral formulas; wherefore, for this, we may
take an arbitrary denominator $a+bx^\theta$ of this differential, so that
we may  have:
\[
ds=\frac{(1-x^\theta)^{\lambda-1}}{a+bx^\theta}(nax^{n-1}dx-
(a(n+\lambda \theta)-bn)x^{n+\theta-1}dx-b(n+\lambda \theta)x^{n+2\theta-1}dx).
\]
Now therefore, by putting after integration $x=1$, we deduce this reduction:
\begin{align*}
na\int \frac{x^{n-1}dx(1-x^\theta)^{\lambda-1}}{a+bx^\theta}=(a(n+\lambda \theta)-bn)
\int \frac{x^{n+\theta-1}dx(1-x^\theta)^{\lambda-1}}{a+bx^\theta}\\
+b(n+\lambda \theta)\int \frac{x^{n+2\theta-1}dx(1-x^\theta)^{\lambda-1}}{a+bx^\theta}.
\end{align*}

33. Here again it is clear what values should arise for $n$ by differentiating
by $\theta$. Moreover, first the value of it should be set $n=\alpha$,
and 
\[
A=\int \frac{x^{\alpha-1}dx(1-x^\theta)^{\lambda-1}}{a+bx^\theta} \,
\textrm{ and } \, B=\int \frac{x^{\alpha+\theta-1}dx(1-x^\theta)^{\lambda-1}}{a+bx^\theta},
\] 
where namely after integration it shall be set $x=1$.
With this having been determined, it may then be made
\begin{align*}
f=\alpha a, f'=(\alpha+\theta)a, \, f''=(\alpha+2\theta)a, \, f'''=(\alpha+3\theta)a, \textrm{ etc.}\\
g=(\alpha+\lambda \theta)a-\alpha b, \, g'=(\alpha+\theta+\lambda \theta)a-(\alpha+\theta)b,\\
g''=(\alpha+2\theta+\lambda \theta)a-(\alpha+2\theta)b, \textrm{ etc.}\\
h=(\alpha+\lambda \theta)b, \, h'=(\alpha+\theta+\lambda \theta)b, \,
h''=(\alpha+\theta+2\lambda \theta)b, \textrm{ etc. [sic]}
\end{align*}
from which will be formed the following continued fraction:
\[
\frac{\alpha aA}{B}=(\alpha+\lambda \theta)a-\alpha b+
\frac{(\alpha+\theta)(\alpha+\lambda \theta)ab}{(\alpha+\theta+
\lambda \theta)a-(\alpha+\theta)b+
\frac{(\alpha+2\theta)(\alpha+\theta+\lambda \theta)ab}{(\alpha+2\theta+\lambda \theta)a-
(\alpha+2\theta)b+(\alpha+3\theta)(\alpha+2\theta+\lambda \theta)ab \textrm{ etc.}}},
\]
the copious expansion of which form we refrain from.

\begin{center}
{\Large VII. The expansion of the formula\\
$s=x^n(e^{\alpha x}(1-x)^\lambda)$}
\end{center}

34. Here therefore it would be
\[
ds=(1-x)^{\lambda-1}(nx^{n-1}dx-(n+\lambda-\alpha)x^ndx-\alpha x^ndx);
\]
then consequently if after integration it is set everywhere $x=1$,
of course in which case it would be $s=1$, we will have this reduction:
\[
n\int x^{n-1}dxe^{\alpha x}(1-x)^{\lambda-1}=(n+\lambda-\alpha)\int x^n dxe^{\alpha x}(1-x)^{\lambda-1}+\alpha \int x^{n+1}dxe^{\alpha x}(1-x)^{\lambda-1}.
\]

35. In this formula therefore, for the exponent $n$ all the values
ascending from unity ought to be taken, where indeed for the minimum value
of this we shall take $n=\delta$, and then the values of the letters
$A$ and $B$ will be able to be elicited from this formula, 
by it being put $x=1$ after integration,
\[
A=\int x^{\delta-1}dxe^{\alpha x}(1-x)^{\lambda-1},
\qquad B=\int x^\delta dxe^{\alpha x}(1-x)^{\lambda-1},
\]
then indeed from these values:
\begin{align*}
f=\delta, \, f'=\delta+1, \, f''=\delta+2, \, f'''=\delta+3, \textrm{ etc.}\\
g=\delta+\lambda-\alpha, \, g'=\delta+1+\lambda-\alpha, \, g''=\delta+2+\lambda-\alpha, \textrm{ etc.}\\
h=\alpha, \, h'=\alpha, \, h''=\alpha, \textrm{ etc.}
\end{align*}
such a continued fraction follows:
\[
\frac{\delta A}{B}=\delta+\lambda-\alpha+\frac{(\delta+1)\alpha}{\delta+1+\lambda-\alpha+
\frac{(\delta+2)\alpha}{\delta+2+\lambda-\alpha+
\frac{(\delta+3)\alpha}{\delta+3+\lambda-\alpha+\textrm{ etc.}}}}
\]
Here it should be noted in particular that the exponents $\lambda$ and
$\delta$ should from necessity be taken as no greater than nothing, 
since otherwise the principal formula $x^ne^{\alpha x}(1-x)^\lambda$
will not vanish in those cases when $x=1$.

36. If a value equal to $1$ were taken for the letters
$\delta$ and $\lambda$,
the case treated above will come forth; and if
integral numbers are assigned to these letters, in the same way
continued fractions will arise,
which by certain operations will be permitted to be reduced to the prior ones.
Truly if either for either one or both of these letters $\delta$ and
$lambda$ fractions are assigned, then forms shall arise which are irreducible
to the prior ones, of which the value would be able to be expressed
in no other way than altogether transcendental quantities. For
if it were $\delta=\frac{1}{2}$ and $\lambda=\frac{1}{2}$, the value
of the letter $A$ will be bound to be obtained from this integral
formula: $A=\frac{e^{\alpha x}dx}{\surd(x-xx)}$ [sic], the integration
of which leads to altogether transcendental quantities, thus so that
the value of such continued fractions comes forth as most abstruse.

\end{document}